\documentclass{proc-l}

\usepackage{amsfonts, amsmath, amssymb, amscd, amsthm, graphicx}

\newtheorem{thm}{Theorem}[section]
\newtheorem{cor}[thm]{Corollary}
\newtheorem{lem}[thm]{Lemma}

\newtheorem{prob}[thm]{Problem}

\theoremstyle{definition}

\theoremstyle{remark}
\newtheorem{rem}[thm]{Remark}

\newcommand{\e }{\varepsilon }

\renewcommand{\ll }{\langle\hspace{-.7mm}\langle }
\newcommand{\rr }{\rangle\hspace{-.7mm}\rangle }
\renewcommand{\pi }{\Pi}
\renewcommand{\sigma }{\Sigma}
\renewcommand{\gamma }{\delta}
\begin{document}

\title{Large groups and their periodic quotients}

\author{A.Yu. Olshanskii}
\address{Department of Mathematics, Vanderbilt University, Nashville, TN
37240, USA} \curraddr{Department of Mathematics, Moscow State
University, Moscow, 119899, Russia}
\email{alexander.olshanskiy@vanderbilt.edu}
\thanks{}

\author{D.V. Osin}
\address{Department of Mathematics, The City College of New York, New York,
NY 10031} \curraddr{} \email{denis.osin@gmail.com}
\thanks{The first author was supported in part by the NSF grants DMS 0245600 and DMS
0455881. The second author was supported in part by the NSF grant
DMS 0605093. Both authors were supported in part by the Russian
Fund for Basic Research grant 05-01-00895.}


\date{}

\maketitle

\begin{abstract}
We first give a short group theoretic proof of the following
result of Lackenby. If $G$ is a large group, $H$ is a finite index
subgroup of $G$ admitting an epimorphism onto a non--cyclic free
group, and $g$ is an element of $H$, then the quotient of $G$ by
the normal subgroup generated by $g^n$ is large for all but
finitely many $n\in \mathbb Z$. In the second part of this note we
use similar methods to show that for every infinite sequence of
primes $(p_1, p_2, \ldots )$, there exists an infinite finitely
generated periodic group $Q$ with descending normal series
$Q=Q_0\rhd Q_1\rhd \ldots $, such that $\bigcap_i Q_i=\{ 1\} $ and
$Q_{i-1}/Q_i$ is either trivial or abelian of exponent $p_i$.
\end{abstract}


\section{Introduction}


Recall that a group $G$ is {\it large} if some finite index
subgroup of $G$ admits a surjective homomorphism onto a
non--cyclic free group. In fact, it is easy to show that if $G$ is
large, then some finite index normal subgroup of $G$ has a free
non--cyclic quotient. Given a subset $S$ of a group $G$, we denote
by $\ll S\rr ^G$ the normal closure of $S$ in $G$. In the recent
paper \cite{Lac}, Lackenby observed that adding higher powered
relations preserves the largeness of a finitely generated group.
More precisely, groups Lackenby proved the
following (under the additional assumption that $G$ was finitely generated).

\begin{thm}\label{main}
Let $G$ be a large group, $H$ a normal subgroup of $G$ of finite
index admitting a surjective homomorphism onto a non--cyclic free
group, $g_1, \ldots , g_k$ elements of $H$. Then the quotient
group $G/\ll g_1^{n}, \ldots, g_k^n \rr ^G$ is large for all but
finitely many $n\in \mathbb N$.
\end{thm}

As noticed in \cite{Lac}, this theorem has interesting
applications to Dehn surgery on $3$--manifolds. It is also
interested from the algebraic point of view since its iterated
applications allow one to construct infinite finitely generated
periodic groups. The proof of Theorem \ref{main} suggested in
\cite{Lac} essentially uses the deep theory related to property
($\tau $) and homology growth developed by Lackenby in his recent
papers (see \cite{Lac} and references therein). In the first part
of this note we provide a short group theoretic proof of Theorem
\ref{main}. In the second part we apply our method to prove a
modified version of Theorem \ref{main} (see Lemma \ref{variant})
and use this version to construct infinite finitely generated
periodic groups with certain specific properties.

More precisely, let $\pi =(p_i)$ be a (finite or infinite)
sequence of primes. In this paper we say that a group $Q$ is a
{\it $\Pi $--graded group} if $Q$ admits a normal series
\begin{equation}\label{ser}
Q=Q_0\rhd Q_1\rhd \ldots ,
\end{equation}
such that $\bigcap\limits_i Q_i=\{ 1\} $ and $Q_{i-1}/Q_i$ is
either trivial or abelian of exponent $p_i$. Note that all
subgroups $Q_i$ have finite index in $Q$ whenever $Q$ is finitely
generated.

\begin{thm}\label{cor}
For any infinite sequence of primes $\pi=(p_i)$, there exists a
finitely generated infinite periodic $\pi $--group.
\end{thm}

Observe that every finitely generated $\pi $--group is residually
finite and if $Q$ is a periodic $\pi $--group, then every element
$x\in Q$ has order $p_{j_1}\cdots p_{j_k}$, where $p_{j_1}, \ldots
,p_{j_m}$ is a subsequence of $\pi $ depending on $x$ (see Lemma
\ref{properties}). For instance, if $\pi =(p,p, \ldots )$, we
obtain an infinite finitely generated periodic residually finite
$p$--group. Such groups were first constructed by Golod
\cite{Golod} and then by Aleshin, Grigorchuk, Gupta - Sidki,
Sushchansky, and others. Note that all these groups are residually
nilpotent.

On the other hand, Theorem \ref{cor} can be applied to find
examples of a quite different nature. Recall that $S$ is a {\it
section} of a group $Q$ if $S$ is a quotient group of a subgroup
of $Q$. Recall also that an infinite group is called {\it just
infinite} if all its proper quotients are finite.

\begin{cor}\label{cor1}
There exists a finitely generated periodic just infinite group $Q$
such that:

\begin{enumerate}
\item Orders of elements of $Q$ are square free.

\item Every section of $Q$ is residually finite. In particular,
$Q$ is residually finite.

\item If $S$ is a finite section of $Q$, then $S$ is solvable and
all Sylow subgroups of $S$ are abelian. In particular, every
nilpotent section of $Q$ is abelian.
\end{enumerate}
\end{cor}

We note that the first examples of finitely generated residually
finite just infinite periodic groups were constructed in
\cite{Grig}. However, groups from \cite{Grig} were $p$--groups. In
particular they were residually nilpotent as well as Golod groups
\cite{Golod} (compare with (c)). To the best of our knowledge no
examples of finitely generated residually finite periodic groups
satisfying either of the properties (a)--(c) were known until now.


\section{Proof of Theorem \ref{main}}


Throughout this paper $F$ denotes a non--abelian free group. We
write $H\twoheadrightarrow F$ if a group $H$ admits an epimorphism
onto some $F$. The main ingredient of our proof is the following.

\begin{thm}[Baumslag, Pride, \cite{BP}] \label{BP}
Suppose that a group $G$ admits a presentation with $n\ge 2$
generators and at most $n-2$ relations. Then $G$ is large.
\end{thm}

\begin{rem}\label{rem}
In fact, Baumslag and Pride proved even a stronger result. Under
the assumptions of Theorem \ref{BP}, they showed that for every
sufficiently large $m\in \mathbb N$ there is a normal subgroup
$H=H(m)\le G$ such that $G/H\cong \mathbb Z/m\mathbb Z$ and
$H\twoheadrightarrow F$ (see \cite{BP}).
\end{rem}

Given elements $g,t$ of a certain group, we use the notation $g^t$
to express $t^{-1}gt$ and $C_G(g)$ to denote the centralizer of $g$
in $G$. The proof of Theorem \ref{main} is based on the following
auxiliary results.

\begin{lem}\label{0}
Let $G$ be a group, $N$ a normal subgroup of $G$, $g\in N$. Let
$T$ denote the set of representatives of the right cosets of $CN$
in $G$, where $C$ is an arbitrary subgroup of $C_G(g)$. Then $\ll
g\rr ^G= \ll Z\rr ^N$, where $Z=\{ g^{t} \, |\, t\in T\}$.
\end{lem}

\begin{proof}
It suffices to show that $g^s\in \ll Z\rr ^N$ for every $s\in G$.
To this end we note that $s=fht$ for some $f \in C\le C_G(g)$,
$h\in N$, and $t\in T$. Hence
$$
g^s=g^{fht}=g^{ht}=g^{th^\prime },
$$
where $h^\prime =h^{t}\in N$. Thus $g^s\in \ll Z\rr ^N$.
\end{proof}

\begin{lem}\label{fi}
For any finite collection of non--trivial elements $g_1, \ldots ,
g_k$ of a free group $F$ and any number $m\in \mathbb N$, there
exists $M\in \mathbb N$ with the following property. For every
$q\ge M$, there is a finite index normal subgroup $N\lhd F$ such
that for all $1\le i\le k$, $g_i^s\notin N$ whenever $1\le s\le
m$, but $g_i^q\in N$.
\end{lem}

\begin{proof}
It suffices to prove the lemma in the case when $F$ is finitely
generated. Let $a_1, ..., a_r$ be a basis of $F$, $A(\mathbb Z,r)$
(respectively $A(\mathbb F_p,r)$) the algebra of formal power
series in non--commutative variables $x_1, \ldots , x_r$ over
$\mathbb Z$ (respectively $\mathbb F_p$), $X$ (respectively $X_p$)
the ideal of $A(\mathbb Z,r)$ (respectively $A(\mathbb F_p,r)$)
generated by $x_1, ..., x_r$. It is well known that the subset
$1+X$ of $A(\mathbb Z,r)$ forms a group with respect to
multiplication, and the map $a_i\to 1+x_i$, $1\le i\le r$, extends
to an embedding $F\to 1+X$ \cite{MKS}.

Clearly there are $l, M_0\in\mathbb N$ such that for any prime
$p\ge M_0$ the natural image of the set $S=\{ g_i^s\, |\, 1\le
i\le k, 1\le s\le m\} $ in  $A(\mathbb F_p,r)/X_p^l$ does not
contain $1$. Without loss of generality we may assume $M_0\ge l$.
Recall that $F$ is residually finite $p$--group for any prime $p$.
Hence for each prime $2\le p\le M_0$, there is a subgroup $N_p$ of
index $p^{j(p)}$ in $F$ for a certain $j(p)\in \mathbb N$ such
that $S\cap N=\emptyset $. Set $M=\prod_{2\le p\le M_0} p^{j(p)}$.
Each number $q\ge M$ is divisible by either $p^{j(p)}$ for some
prime $2\le p\le M_0$, or by a prime $p> M_0\ge l$. In the first
case, we set $N=N_p$. In the second case the inequality $p>l$
implies that the image of every element of $1+X$ in $A(\mathbb
F_p,r)/X_p^l$ has order $p$. Therefore, we can we set $N$ to be
the kernel of the homomorphism of $F$ onto its image in $A(\mathbb
F_p,r)/X_p^l$.
\end{proof}

\begin{lem}\label{1}
Let $F$ be a free group of rank $r\ge 2$, $g_1, \ldots g_k$
arbitrary elements of $F$. Then  $\overline{F}=F/\ll g_1^{q}, \ldots
, g_k^{q}\rr ^{F}$ is large for all but finitely many $q\in \mathbb
N$.
\end{lem}

\begin{proof}
Without loss of generality we may assume that $g_1, \ldots , g_k$
are nontrivial. By Lemma \ref{fi}, there exists $M\in \mathbb N$
such that for any $q\ge M$ there is a finite index normal subgroup
$N\lhd F$ such that for all $i=1, \ldots , k$, $g_i^s\notin N$ if
$s=1, \ldots , k$, but $g_i^q\in N$. Note that the image of $g_i$
in $F/N$ has order at least $k+1$. Therefore the index of $\langle
g_i\rangle N$ in $F$ is at most $j/(k+1)$, where $j=|G:N|$.

According to Lemma \ref{0}, the image $\overline{N}$ of the
subgroup $N$ in $\overline{F}$ is isomorphic to $N/\ll Z\rr ^N$,
where $Z=\bigcup\limits_{i=1}^k \{ (g_i^{q})^{t} \, |\, t\in
T_i\}$ and $T_i$ is the set of representatives of the right cosets
of $\langle g_i\rangle N$ in $F$. Therefore, $\overline{N}$ admits
a presentation with
$$
rank\, N= 1+(r-1)j\ge 1+j
$$
generators and
$$
\# Z=\sum\limits_{i=1}^k \# T_i =\sum\limits_{i=1}^k |F:\langle
g_i\rangle N |\le \frac{kj}{k+1}<j
$$
relations. Hence $\overline{N}$ is large by Theorem \ref{BP}. As
$R$ is of finite index in $\overline{F}$, the group $\overline{F}$
is also large.
\end{proof}

\begin{proof}[Proof of Theorem \ref{main}]
Let $H$ be a normal subgroup of finite index in $G$ such that
admitting a homomorphism $\e\colon H\to F$ onto a non--cyclic free
group. By Lemma \ref{0}, for every $n\in \mathbb N$, the image
$\overline{H}$ of $H$ in $G/\ll g_1^{n}, \ldots , g_k^n\rr ^G$ is
isomorphic to $H/\ll Z^n\rr ^H$, where $Z$ is some finite subset
of $H$ consisting of conjugates of $g_i$'s and $Z^n=\{z^n | z\in
Z\}$ . Thus $\overline{H}$ admits a surjective homomorphism onto
$F/\ll \e(Z^n) \rr ^{F}=F/\ll \e(Z)^n \rr ^{F}$, which is large
for all but finitely many $n$ according to Lemma \ref{1} applied
to the set $\e(Z)$. Hence so are $\overline{H}$ and $G/\ll
g^{n}\rr ^G$.
\end{proof}


\section{Constructing periodic $\pi $--graded groups}


To each group $G$ and each sequence of primes $\sigma =(q_i)$, we
associate a sequence of characteristic subgroups $\gamma_i^\sigma
(G)$ of $G$ defined by $\gamma _0^\sigma (G)=G$ and
$$\gamma _i^\sigma (G)=[\gamma _{i-1}^\sigma (G), \gamma
_{i-1}^\sigma (G)] \left(\gamma _{i-1}^\sigma (G)\right)^{q_i}.$$
It is easy to prove by induction that these subgroups have finite
index in $G$ whenever $G$ is finitely generated. We need an
auxiliary lemma, which follows immediately from a result of Levi
(see \cite[Lemma 21.61]{HN}).

\begin{lem}\label{Levi}
For any finite subset $S$ of nontrivial elements of a free group $F$, there
exists $D\in \mathbb N$ such that for any infinite sequence of primes $\sigma $, we have
$\gamma _d^\sigma (F)\cap S=\emptyset $ for all $d\ge D$.
\end{lem}

Throughout the rest of this section we fix an arbitrary infinite
sequence $\pi =(p_i)$ of primes and denote by $\omega _k$ the
subsequence $(p_{k+1}, p_{k+2}, \ldots )$ of $\pi $.

\begin{lem}\label{variant}
Let $G$ be a finitely generated group. Suppose that $\gamma _r^\pi
(G)\twoheadrightarrow F$ for some $r$. Then for any element $g\in
\gamma _r^\pi (G)$ there is $m\ge r$ such that if $g^n\in \gamma
_m^\pi (G)$, then $\gamma _m^\pi (G)/\ll g^{n}\rr
^G\twoheadrightarrow F$.
\end{lem}

\begin{proof}
By Lemma \ref{0}, $\gamma _r^\pi (G)/\ll g^{n}\rr ^G$ is
isomorphic to $\gamma _r^\pi (G)/\ll h_1^{n}, \ldots , h_k^{n}\rr
^{\gamma _r^\pi (G)}$, where $h_1, \ldots , h_k$ are some
conjugates of $g$. Thus $\gamma _r^\pi (G)/\ll g^{n}\rr ^G$
surjects onto $\overline{F}=F/\ll g_1^{n}, \ldots , g_k^{n} \rr
^{F}$, where $g_i$ stands for the image of $h_i$ in $F$. Without
loss of generality we may assume that $g_1, \ldots , g_k$ are
non--trivial. By Lemma \ref{Levi}, there exists $d$ such that
$g_i^s\notin \gamma _d^{\omega _r} (F)$ for $1\le i, s\le k$. We
denote the subgroup $\gamma _d^{\omega _r}(F)$ by $N$.

Assume that $g^n\in \gamma _{r+d}^\pi (G)$. Then $g_i^n\in N$ for
all $i$. We now repeat word--for--word the arguments from the
proof of Lemma \ref{1} and conclude that the image $\overline{N}$
of $N$ in $\overline{F}$ admits a presentation with $2$ more
generators than relations. According to Remark \ref{rem}, there is
a subgroup $M\lhd \overline{N}$ such that $\overline{N}/M$ is
cyclic of order $p_{r+d+1}\cdots p_{r+d+c}$ for some $c\ge 1$ and
$M\twoheadrightarrow F$. Clearly $\gamma _{c+d}^{\omega _{r}}
(\overline{N})$ is a finite index subgroup of $M$ and hence
$\gamma _{c+d}^{\omega _{r}}(\overline{N})\twoheadrightarrow F$.
It remains to note that $\gamma _{r+d+c}^\pi (G)/ \ll g^{n}\rr ^G$
surjects onto $\gamma _{c+d}^{\omega _{r}}(\overline{N})$ and set
$m=r+d+c$.
\end{proof}

\begin{proof}[Proof of Theorem  \ref{cor}]
We enumerate all elements of the free group $F=\{f_1, f_2, \ldots
\} $ and construct a sequence of quotients $G_i$ of $F$ as
follows. Let $G_0=F$ and suppose that we have already constructed
$G_i$ such that
\begin{equation}\label{gri}
\gamma _{r_i}^\pi (G_i)\twoheadrightarrow F
\end{equation}
for some $r_i\in \mathbb N$. Let also $g_{i+1}$ denote the image
of $(f_{i+1})^{p_1p_2\cdots p_{r_i}}$ in $G_i$. Then $g_{i+1}\in
\delta_{r_i}^\Pi (G)$ and we can choose $r_{i+1}> r_i$ and $n
=n(i)\in \mathbb N$ according to Lemma \ref{variant} so that
$g_{i+1}^n\in \gamma _{r_{i+1}}^\pi (G_i)$ and $\gamma
_{r_{i+1}}^\pi (G_i)/ \ll g_{i+1}^n\rr ^{G_i}\twoheadrightarrow
F$.  We set $G_{i+1}=G_i/\ll g_{i+1}^n\rr ^{G_i}$. Clearly $\gamma
_{r_{i+1}} ^\pi (G_{i+1})\twoheadrightarrow F$. Observe also that
(\ref{gri}) and inequality $r_{i+1}>r_i$ imply that $\gamma
_{r_{i+1}}^\pi (G_i)$ is a proper subgroup of $\gamma _{r_{i}}^\pi
(G_i)$. Therefore,
\begin{equation}\label{order}
|G_{i+1}/\gamma _{r_{i+1}}^\pi (G_{i+1})|=|G_i/\gamma
_{r_{i+1}}^\pi (G_i)|>|G_i/\gamma _{r_{i}}^\pi (G_i)|.
\end{equation}

Let $K_i$ denote the kernel of the natural homomorphism $F\to
G_i$. Clearly $G=F/\bigcup\limits_{i=1}^\infty K_i $ is a periodic
group. Further set $Q=G/\bigcap\limits_{i=1}^\infty \gamma ^\pi
_{r_i}(G)$. Obviously $Q$ is a $\pi $--group. To show that $Q$ is
infinite we observe that $Ker (G_i\to G)\le \gamma _{r_i}^\pi
(G_i)$ for every $i$. Hence $$G/\gamma _{r_i}^\pi (G)\cong
G_i/\gamma ^\pi _{r_i}(G_i).$$ Inequality (\ref{order}) now
implies $|G/\gamma ^\pi _{r_i}(G)|\to \infty $ as $i\to \infty $.
Therefore, the group $Q$ is infinite as it surjects onto $G/\gamma
^\pi _{r_i}(G)$ for every $i$.
\end{proof}

To derive Corollary \ref{cor1} we first list some general
properties of $\pi $--graded groups in the case when $\pi $
consists of distinct primes.

\begin{lem}\label{properties}
Suppose that all primes in the sequence $\pi =(p_1,p_2, \ldots )$
are distinct and $Q$ is a periodic $\pi $--graded group with
normal series (\ref{ser}). Then the following holds.
\begin{enumerate}
\item[(a)] Every element $x\in Q$ has order $p_{j_1} \ldots
p_{j_k}$, where $p_{j_1}, \ldots , p_{j_k}$ is a subsequence of
$\pi $. Moreover, if $x\in Q_i$, then $j_1>i$.

\item[(b)] Every subgroup of $Q$ is a $\pi $--graded group.

\item[(c)] Every quotient $\overline{Q}$ of $Q$ is a $\pi
$--graded group.

\item[(d)] If $Q$ is finitely generated, then every section of $Q$
is residually finite.

\item[(e)] Every finite section of $Q$ is solvable. Every
nilpotent section of $Q$ is abelian.
\end{enumerate}
\end{lem}

\begin{proof}
To prove (a) we note that for every element $x\in Q$, there exists
$m$ such that $Q_m\cap \langle x\rangle =\{ 1\} $. Hence the order
of $x$ in $Q$ coincides with the order of the image of $x$ in the
$\pi $--graded group $Q/Q_m$. Thus it suffices to prove (a) for
$\pi $--graded groups admitting finite series of type (\ref{ser}).
This can easily be done by induction on the length of the series.

Assertion (b) is trivial. To prove (c) it suffices to note that
the intersection $I$ of the images of subgroups $Q_i$ in
$\overline{Q}$ is periodic, but $I$ can not have elements of any
prime order according to (a).

If $R$ is a subgroup of $Q$ and $\overline{R}$ is a quotient of
$R$, then the same argument as above shows that the intersection
of images of $\overline{R}_i=R\cap Q_i$ in $\overline{R}$ is
trivial. Since $Q$ is finitely generated,
$\overline{R}/\overline{R}_i\cong RQ_i/Q_i\le Q/Q_i$ is finite.
This proves (c).

Note that every section of $Q$ is a $\pi $--graded group by (b)
and (c). Every finite $\pi$--graded group admits a series
(\ref{ser}) of finite length, hence every such a group is
solvable. Finally we recall that periodic nilpotent groups are
locally finite and finite nilpotent groups are direct products of
their Sylow subgroups. Obviously every finite $\pi $--graded
$p$--group is abelian whenever $\pi $ consists of distinct primes.
This implies (d).
\end{proof}

\begin{proof}[Proof of Corollary \ref{cor1}]
Suppose that $\pi $ consists of distinct primes and $Q$ is a $\pi
$--graded group provided by Theorem \ref{cor}. As $Q$ is finitely
generated and infinite, by the Zorn Lemma there is a just infinite
quotient $\overline{Q}$ of $Q$. The group $\overline{Q}$ is also a
$\pi $--group and have all required properties by Lemma
\ref{properties}.
\end{proof}

Finally we mention one question motivated by Theorem \ref{cor} and
Corollary \ref{cor1}.

\begin{prob}\label{prob}
Does there exist an infinite finitely generated periodic
residually (finite simple) group?
\end{prob}

Another interesting question is whether there exists an infinite
finitely generated residually finite group that is also residually
simple. Problem \ref{prob} is related to one of the main open
questions in the theory of hyperbolic groups, that is, whether
every hyperbolic groups is residually finite. If this question has
positive answer, our method can be used to construct finitely
generated infinite periodic residually (finite simple) groups.
Here is the sketch of the proof.

In \cite{Ols00}, the first author observed that if all hyperbolic
groups are residually finite, then every non--elementary
hyperbolic group has infinitely many finite simple quotients.
Recall also that if $G$ is a non--elementary hyperbolic group and
$g\in G$, then $G/\ll g^n\rr ^G$ is also non--elementary
hyperbolic provided $n$ is big enough \cite{Ols93}. Starting with
a finitely generated non--abelian free group $F=\{ f_1, f_2,
\ldots \} $ and assuming that every hyperbolic group is residually
finite, we can construct a sequence of non--elementary hyperbolic
quotients $G_i$ of $F$ and subgroups $N_i\lhd G_i$ such that
$G_{i+1}=G_i/\ll g_{i+1}^n\rr ^{G_i}$ for some $n=n(i)$, $G_i/N_i$
is finite simple, $|G_i/N_i|\to \infty $ when $i\to \infty $, and
for every $i$, $g_{i+1}^n$ belongs to the intersection of the
images of all subgroups $N_0,N_1, \ldots , N_i$ in $G_i$. Here
$g_i$ is the image of $f_i$ in $G_i$ as above. Then we define the
quotient group $G$ of $F$ as in the proof of Theorem \ref{cor} and
denote by $Q$ the quotient of $G$ by the intersection of images of
all subgroups $N_i$, $i\in \mathbb N\cup \{ 0\} $ in $G$. Clearly
$Q$ is finitely generated, periodic, residually (finite simple),
and infinite for the same reason as above.

\end{document}